\documentclass[11pt]{article}
\usepackage{amsmath, epsfig}
\usepackage{mathrsfs}
\usepackage{latexsym}
\usepackage{amssymb}

\makeatother

\usepackage{graphicx}

\makeatother

\begin{document}

\title{Further results on computation of topological indices of certain networks}

\footnotetext{This work is partially supported by National Natural Science Foundation of China (nos. 11471016,
11401004, 11571134, 11371162), Anhui Provincial Natural Science Foundation (nos. KJ2015A331,
KJ2013B105, 1408085QA03), the Self-determined Research Funds of Central China Normal
University from the colleges basic research and operation of MOE, and the Summer Graduate Research Assistantship Program of Graduate School at the University of Mississippi. }

\author{Shaohui Wang$^a$,  Jia-Bao Liu$^b$\footnote{  Corresponding author.  \newline E-mail addresses: S. Wang  (swang4@go.olemiss.edu, shaohuiwang@yahoo.com), J.B.  Liu (liujiabaoad@163.com), C. Wang (wcxiang@mailccnu.edu.cn), S. Hayat(sakander1566@gmail.com)}, $~$Chunxiang Wang$^c$, Sakander Hayat$^d$\\
\small\emph {a. Department of Mathematics, The University of Mississippi, University, MS 38677, USA}
\\
\small\emph {b. School of Mathematics and Physics, Anhui Jianzhu University, Hefei, 230601, P.R. China}\\
\small\emph {c. School of Mathematics and Statistics, Central China Normal University, Wuhan, 430079,  P.R. China}\\
\small\emph {d. School of Mathematical Sciences, University of Science and Technology of China,
Hefei, 230026, China}}
\date{}
\maketitle

\begin{abstract}
There are various topological indices such as degree based
topological indices, distance based topological indices and
counting related topological indices etc. These topological
indices correlate certain physicochemical properties such as
boiling point, stability of chemical compounds. In this paper, we compute the sum-connectivity index and multiplicative Zagreb indices for certain networks of chemical importance like silicate networks, hexagonal networks, oxide networks, and honeycomb networks.
Moreover, a comparative study using computer-based graphs has been made to clarify their nature for these families of networks.

\end{abstract}

\section{Introduction}
 Throughout this paper $G = (V, E)$ is a finite, simple and connected graph with vertex set $ V(G)$ and edge set $E(G)$.
For other undefined notations, readers may refer
to~\cite{Bondy1976,LiuP2015,Liu2015}. Chemical graph theory is a
branch of graph theory whose focus of interest is finding
topological indices of chemical graphs which correlate well with
chemical properties of the chemical molecules. A topological index
is a numerical parameter mathematically derived from the graph
structure.  The topological indices have been found to be useful for
establishing correlations between the structure of a molecular
compound and its physicochemical properties or biological
activity~\cite{dneg2016,dengh2015, Khadikar2000}.

Let $d(u),  ~d(v)$ be the degrees of the vertices $u$ and $v$ respectively. One
of the topological indices used in mathematical chemistry is that of
the so-called degree-based topological indices, which are defined
in terms of the degrees of the vertices of a graph. For instance,
the first and second Zagreb indices of $G$ are respectively
defined as
 $$ M_1(G) = \sum_{u \in V(G)} d(u)^2,
~M_2(G) = \sum_{uv \in E(G)} d(u)d(v).$$
The background and applications of Zagreb indices can be found in \cite{BRu1975,Gutman2015,GutmanN,Siddiqui2016}.

 In the
1980s, Narumi and Katayama~\cite{Nar} characterized  the
structural isomers of saturated hydrocarbons and considered the
product $ NK(G) = \prod_{v \in V(G)} d(v), $ which is called the
NK index.
Two fairly new indices with higher prediction
ability~\cite{RT2010}, named the first and second multiplicative
Zagreb indices~\cite{Gutman2011}, are respectively defined as
 $$\Pi_{1,c} (G) = \prod_{v \in V(G)}d(v)^c,  ~\Pi_2(G) = \prod_{uv \in E(G)}
 d(u)d(v).$$
Obviously, the first multiplicative Zagreb index is the power of the NK index.
 Moreover, the second multiplicative Zagreb index can be rewritten as $\prod_2(G) = \prod_{u \in V(G)}d(u)^{d(u)}$. The properties of $\prod_{1,c}(G), \prod_{2}(G)$  for some chemical structures
have been studied extensively in \cite{Wang2015, Eliasi2012, Gutman20015, Shi2015}

The multiplicative version of ordinary first
Zagreb index $M_1$ is defined as
$$\Pi_1^*(G) = \prod_{uv \in E(G)} (d(u)+d(v)).$$
The general sum-connectivity index is defined as
$$\chi_{\alpha}(G) = \sum_{uv \in E(G)} (d(u)+d(v))^{\alpha}, \text{where} ~\alpha~ \text{is a real number,}$$
and given  with the purpose of extending the
classical sum-connectivity index, $\chi_{\frac{1}{2}}(G).$ The properties of $\prod_1^*(G), \chi_{\alpha}(G)$ polynomials for some chemical structures have been studied in \cite{Eliasi2012, Zhou2010}.
 For other work
on topological indices, the readers are referred
to~\cite{Berrocal2016,mr2014,mr2013,mrf2013,Gutman2014, Kazemi2016, liu2015,liuw2015,xxl,Nar, T2016,wang2016,Xu,xxz}.

Silicates are the most important elements of
 Earth's crust, as well as the other terrestrial planets, rocky moons, and asteroids.  Sand, Portland cement, and thousands of minerals are constituted by silicates.
The tetrahedron $(SiO_4)$ is a basic unit of silicates,  in which the central vertex is silicon vertex and the corner vertices are oxygen vertices, see Figure 1. 

One of most interesting topics is to obtain the quantitative structure activity/
property/toxicity relationships \cite{emeric1, emeric, Hayat2014, RAJAN2012, Siddiquia2016}. In this paper, we explore the silicate, chain silicate,
hexagonal, oxide and honeycomb networks and provide exact polynomials of
$\chi_{\alpha}(G),$ $\prod_1^*(G),$   $\prod_{1,c}(G)$
and $\prod_2(G)$ for these networks. Furthermore, we use the pictures to compare
the chemical indices on these silicate-related networks.

\begin{figure}[htbp]
    \centering
    \includegraphics[width=1.3in]{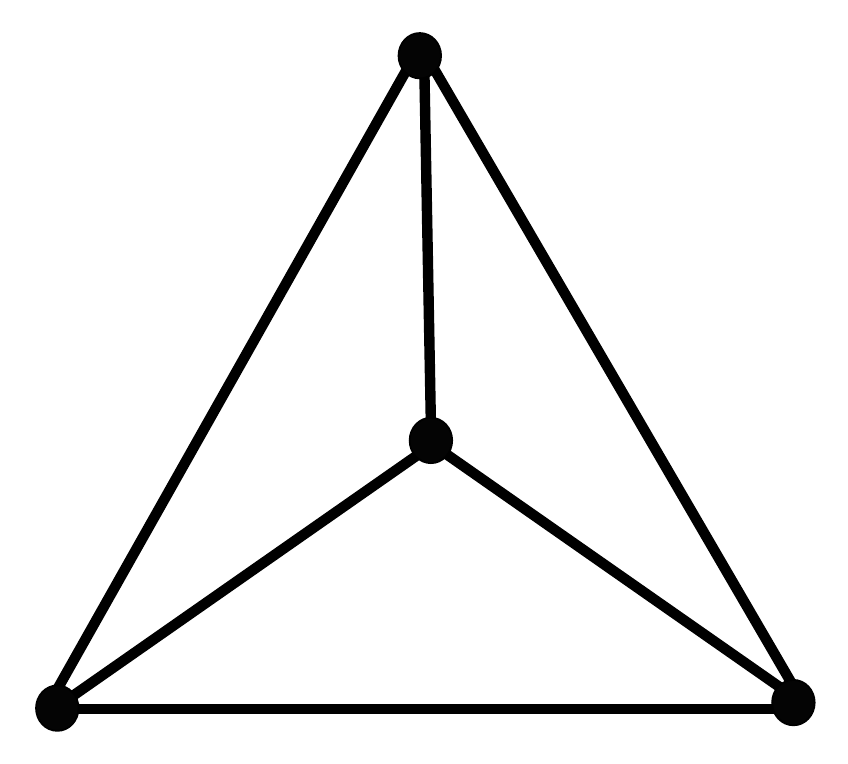}
    \caption{ $SiO_4$ Tetrahedron.}
    \label{fig: te}
\end{figure}

\begin{figure}[htbp]
    \centering
    \includegraphics[width=2.5in]{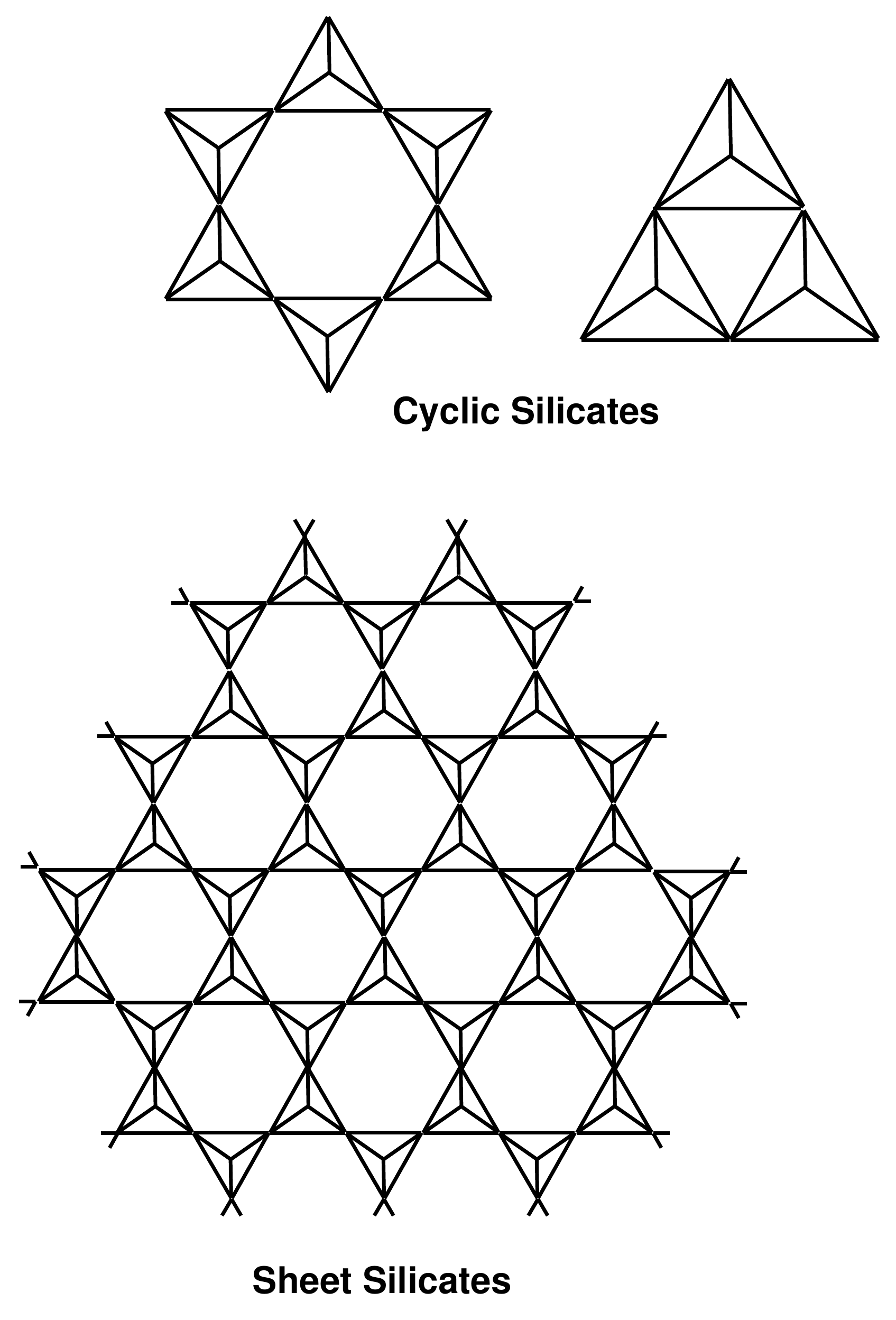}
    \caption{Different sheet silicates.}
    \label{fig:sample_figure}
\end{figure}

\section{Main results and discussion}
In this section, we explore the silicate, chain silicate,
hexagonal, oxide and honeycomb networks and provide exact formulas of
$\chi_{\alpha}(G),$ $\prod_1^*(G),$   $\prod_{1,c}(G)$
and $\prod_2(G)$ for these networks. In order to
compute certain topological indices of thee silicate networks, we will divide the
vertex set and the edge set based on degrees of end vertices of each edge of the graph. In addition, we utilized some techniques and figures(Figs. 1 - 7) in \cite{Hayat2014, RAJAN2012}.

We first give the general remarks how to compute the relation between the numbers of vertices and edges from the literature. In
 general, we follow the approaches of Gutman, Deutsch and Kla\u{v}zar \cite{emeric1,emeric,Gutmans} for which some definitions
are needed first.

\vskip 2mm \noindent  {\bf Lemma 2.0.1. \cite{emeric1,emeric,Gutmans}} \emph{Given that $G$ is a connected graph with $n=|V(G)|$ and $m=|E(G)|$.
Let $n_i$ be the number of vertices of degree $i$ and $m_{i,j}$ be the number of edges with end vertices of degree $i$ and $j$, $i,j \geq 1$. Then
$$\begin{array}{rcl}
&n_1+n_2+n_3+n_4+n_5+n_6  &=n \\
& m_{12}+m_{13}+m_{14}+m_{15}+m_{16} & =n_1,\\
& m_{21}+2m_{22}+m_{23}+m_{24}+m_{25}+m_{26} & =2n_2,\\
& m_{31}+ m_{32}+2m_{33}+m_{34}+m_{35}+m_{36} & =3n_3,\\
& m_{41}+ m_{22}+m_{43}+2m_{44}+m_{45}+m_{46} & =4n_4,\\
& m_{51}+ m_{52}+m_{53}+m_{54}+2m_{55}+m_{56} & =5n_5,\\
& m_{61}+ m_{62}+m_{63}+m_{64}+m_{65}+2m_{66} & =6n_6,\\
& n_1+2n_2+3n_3+4n_4+5n_5+6n_6 & =2m.\\
\end{array}$$}

\subsection{Silicate networks}

 Silicates are the largest, very interesting and most complicated minerals by far.
Silicates are obtained by fusing metal oxides or metal carbonates
with sand~\cite{Hayat2014}, see Figure 2. A silicate
network of dimension $n$ symbolizes as $SL_n$, where $n$ is the number of hexagons between the center and boundary of
$SL_n$. A  silicate network of dimension two is shown in Figure 3. The number of vertices in $SL_n$ are
$15n^2+ 3n$, and the number of edges are $36n^2$.

\begin{figure}[htbp]
    \centering
    \includegraphics[width=2.5in]{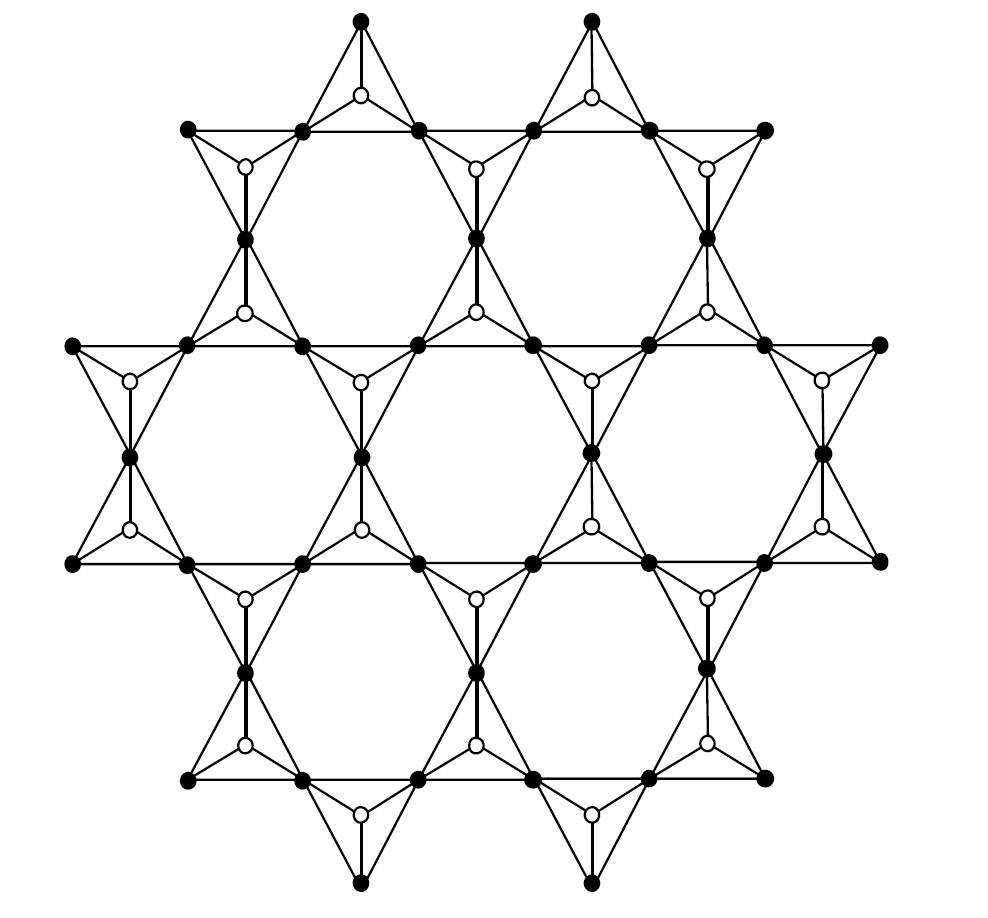}
    \caption{Silicate network of dimension two, the solid vertices denote the oxygen atoms whereas the plain vertices represent the silicon atoms.}
    \label{fig:sample_figure}
\end{figure}

In the following theorem, the exact formulas of
$\chi_{\alpha}(SL_n),$ $\prod_1^*(SL_n),$   $\prod_{1,c}(SL_n)$
and $\prod_2(SL_n)$ for silicate networks are computed.

\vskip 2mm \noindent  {\bf Theorem 2.1.1. } \emph{Consider the
 silicate networks $SL_n$, then the indices of
$\chi_{\alpha}(SL_n),$ $ \prod_1^*(SL_n),$ $ \prod_{1,c}(SL_n)$
and $\prod_2(SL_n)$ are equal to
$$\begin{array}{rcl}
&~\prod_{1,c}(SL_n)&  = 2^{9n^2c -3nc}3^{15n^2c+3nc}, \\
& \prod_2(SL_n)& = 2^{54n^2-18n}3^{72n^2},\\
 &\chi_{\alpha}(SL_n)& = n6^{\alpha+1} +3^{2\alpha+1}(6n^2+2n)+2^{2\alpha+1}3^{\alpha+1}(3n^2-2n), \\
&\prod^*_1(SL_n)& =  2^{36n^2-18n}3^{54n^2+6n}.
\end{array}$$}
 \vskip 2mm \noindent  {\bf Proof. }
Let $G$ be the graph of silicate network $SL_n$ with $| V(SL_n) | = 15
n^2+3n$ and $| E(SL_n) | = 36n^2$. By the properties of $SL_n$,  the vertex partition of silicate networks $SL_n$ based on the
degree  of vertices comprises the vertices of degrees 3 and 6. The vertices on the center of $SiO_4$ are of degree 3,
  the vertices on both the boundary of $SL_n$ and its $SiO_4$ are degree 3 and others are of degree 6.
  Table 1  provides such partition for
$SL_n$ below.

\begin{table}[thb]
\center
\caption{Vertex partition of $SL_n$ based on degrees of vertices.}
\label{tab:CS}
\vspace{0.2in}
\begin{tabular}{|lc  | c|c|c||c ||}
\hline
 \multicolumn{2}{|c|}{Degrees} &\multicolumn{ 1}{c|}{$\;\;\;\;\;$$3$$\;\;\;\;\;$}  &\multicolumn{ 1}{c|}{$\;\;\;\;\;$$6$$\;\;\;\;\;$}  \\
\hline\hline
  & Number of vertices     & \multicolumn{ 1}{c|} {$6n^2+6n$} & \multicolumn{ 1}{c|}{$9n^2-3n$}  \\
\hline
\end{tabular}
\end{table}

By Table 1, $\prod_{1,c}(G) = \prod_{v \in V(G)}d(v)^c$ and  $\prod_2(G) = \prod_{u \in V(G)}
d(u)^{d(u)}$, we have
$$\begin{array}{rcl}
&~\prod_{1,c}(SL_n)&  =3^{c(6n^2+6n)}6^{c(9n^2-3n)}= 2^{9n^2c -3nc}3^{15n^2c+3nc}, \\
& \prod_2(SL_n)& = 3^{3(6n^2+6n)}6^{6(9n^2-3n)} = 2^{54n^2-18n}3^{72n^2}.
\end{array}$$
From the graphic properties of a silicate network, there are three
types of edges based on the degree of the vertices of each edge.
By Table 1, we get $$n_3=6n^2+6n,~~~n_6=9n^2-3n.$$ Note that
$m_{3,3} + m_{3,6} + m_{6,6} = 36n^2 =\mid E(SL_n)\mid.$ In
addition, by Lemma 2.0.1, one can obtain
$$
2m_{3,3} + m_{3,6} = 3n_3,~
 m_{3,6} + 2m_{6,6} = 6n_6.
$$
Consequently,
$$ 
m_{3,3} = 6n,~
m_{3,6} = 18n^2+6n, ~
 m_{6,6} = 18n^2-12n. 
$$
 The following table gives the three types and the number
of edges in each type for $SL_n$.
\begin{table}[thb]
\center
\caption{Edge partition of $SL_n$ based on degrees of end vertices of each edge.}
\label{tab:CS}
\vspace{0.2in}
\begin{tabular}{|lc  | c|c|c||c ||}
\hline
 \multicolumn{2}{|c|}{($d(u),d(v)$) where $u,v \in E(G)$} &\multicolumn{ 1}{c|}{$\;\;\;\;\;$$(3,3)$$\;\;\;\;\;$}
   &\multicolumn{ 1}{c|}{$(3,6)$}&\multicolumn{ 1}{c|}{$(6,6)$}  \\
\hline\hline
  & Number of edges     & \multicolumn{ 1}{c|} {$6n$} & \multicolumn{ 1}{c|}{$\;\;\;\;\;$$18n^2+6n$} &
   \multicolumn{ 1}{c|}{$\;\;\;\;\;$$18n^2-12n$} \\
\hline
\end{tabular}
\end{table}

By Table 2, $\chi_{\alpha}(G) = \sum_{uv \in E(G)} (d(u)+d(v))^{\alpha}$ and $\prod_1^*(G) = \prod_{uv \in E(G)} (d(u)+d(v))$,
 we obtain $$\begin{array}{rcl}
 &\chi_{\alpha}(SL_n)& = 6^{\alpha}\cdot6n+ 9^{\alpha}(18n^2+6n)+12^{\alpha} (18n^2-12n) \\ &&= n6^{\alpha+1} +3^{2\alpha+1}(6n^2+2n)+2^{2\alpha+1}3^{\alpha+1}(3n^2-2n), \\
&\prod^*_1(SL_n)& =  2^{36n^2-18n}3^{54n^2+6n}.
\end{array}$$
This completes the proof. $\hfill\Box$


\subsection{ Chain silicate networks}
Next, we provide another family of silicate networks named as chain
silicate networks and then compute its certain degree based
topological indices.
 Here we provide chain silicate networks $CS_n$ of dimension $n$ as follows:
  A chain silicate network of dimension $n$ symbolizes as $CS_n$ is obtained by arranging $n$ tetrahedra linearly, see Figure 4.
  The number
of vertices in  $CS_n$  are $3n+ 1$ and number of edges are $6n$.

\begin{figure}[htbp]
    \centering
    \includegraphics[width=4.0in]{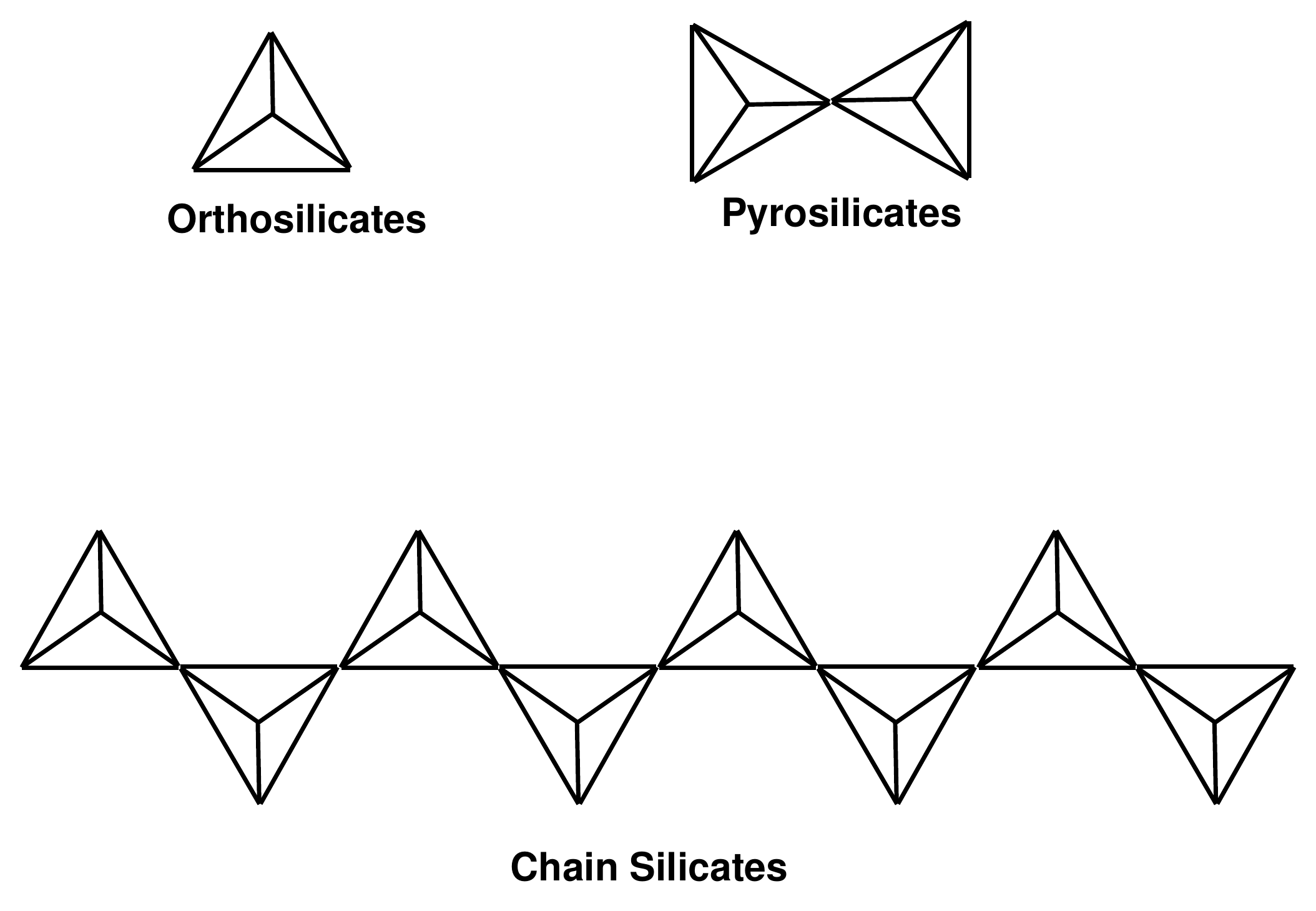}
    \caption{Ortho, Pyro and chain silicates.}
    \label{fig:sample_figure}
\end{figure}

\vskip 2mm \noindent  {\bf Theorem 2.2.1. } \emph{Consider the
chain silicate networks $CS_n$, then  the indices of
$\chi_{\alpha}(CS_n),$ $\prod_1^*(CS_n),$ $ \prod_{1,c}(CS_n)$ and
$\prod_2(CS_n)$ are equal to
$$\begin{array}{rcl}
&~\prod_{1,c}(CS_n)&  = 2^{nc-c}3^{3nc+c}, \\
& \prod_2(CS_n)& = 2^{6n-6}3^{12n},\\
 &\chi_{\alpha}(CS_n)& = \left\{ \begin{array}{rcl}
6^{\alpha +1,}~~~~~~~~~~~~~~~~~~~~~~~~~~~~~~~~~~~~~~~~~~~~~ & \mbox{if} &  n =1, \\
2^{\alpha} 3^{\alpha}(n+4)+  3^{2\alpha} (4n-2)+2^{2\alpha}3^{\alpha}(n-2),  & \mbox{if} &  n \geq 2,
\end{array}\right. \\
&\prod^*_1(CS_n)& =  \left\{ \begin{array}{rcl}
46656,~~~~ ~~~~~~~~~~~~~~~~~~~~~~~~~~~~ & \mbox{if} &  n =1, \\
2^{3n}3^{10n-2},~~~~~~~~~~~~~~~~~~~~~~~~~~~& \mbox{if} & n \geq 2.
\end{array}\right.
\end{array}$$
}

 \vskip 2mm \noindent {\bf Proof. }
Let $G$ be the graph of chain silicate network $CS_n$ with
$| V(CS_n) | = 3n+ 1$ and $| E(CS_n) | = 6n$. Based on  the graphic
properties with the degree  of vertices,  the vertex
partition of silicate networks $CS_n$ are vertices of degree 3 and that of 6.   Table 3  explains such
partition for $CS_n$.

\begin{table}[thb]
\center
\caption{Vertex partition of $CS_n$ based on degrees of vertices.}
\label{tab:CS}
\vspace{0.2in}
\begin{tabular}{|lc  | c|c|c||c ||}
\hline
 \multicolumn{2}{|c|}{Degrees} &\multicolumn{ 1}{c|}{$\;\;\;\;\;$$3$$\;\;\;\;\;$}  &\multicolumn{ 1}{c|}{$\;\;\;\;\;$$6$$\;\;\;\;\;$}  \\
\hline\hline
  & Number of vertices     & \multicolumn{ 1}{c|} {$2n+2$} & \multicolumn{ 1}{c|}{$n-1$}  \\
\hline
\end{tabular}
\end{table}

By Table 3, $\prod_{1,c}(G) = \prod_{v \in V(G)}d(v)^c$ and  $\prod_2(G) = \prod_{u \in V(G)}
d(u)^{d(u)}$, we have
$$\begin{array}{rcl}
&~\prod_{1,c}(CS_n)&  =3^{c(2n+2)}6^{c(n-1)}=  2^{nc-c}3^{3nc+c}, \\
& \prod_2(CS_n)& = 3^{3(2n+2)}6^{6(n-1)}=2^{6n-6}3^{12n}.
\end{array}$$

The edge partition of silicate networks $CS_n$ can be obtained by
the degree sum of vertices for each edge of the graph.
 Note that
$m_{3,3} + m_{3,6} + m_{6,6} = n + 4 + 4n - 2 + n - 2 = 6n =\mid
E(CS_n)\mid.$ In addition, by Lemma 2.0.1, one can obtain
$$ 
 2m_{3,3} + m_{3,6} = 3n_3,~
 m_{3,6} + 2m_{6,6} = 6n_6.
$$
 Combining with $n_3=2n+2, ~~~n_4=n-1,$ we arrive to
 $$m_{3,3}=n+4, ~~m_{3,6}=4n-2, ~~m_{6,6}=n-2.$$
 Table 4 gives such edges partition for $CS_n$.

\begin{table}[thb]
\center
\caption{Edge partition of $CS_n$ based on degrees of end vertices of each edge.}
\label{tab:CS}
\vspace{0.2in}
\begin{tabular}{|lc  | c|c|c|c ||}
\hline
 \multicolumn{2}{|c|}{$d(u),d(v)$) where $u,v \in E(G)$} &\multicolumn{ 1}{c|}{$\;\;\;\;\;$$(3,3)$$\;\;\;\;\;$}  &\multicolumn{ 1}{c|}{$(3,6)$}  &\multicolumn{ 1}{c|}{$(6,6)$}   \\
\hline\hline
      Number of edges  & $n =1 $ &  6  &0 &0\\
 & $n\geq2$ & $n+4$&$4n-2$& $n-2$ \\
\hline
\end{tabular}
\end{table}

By Table 4, $\chi_{\alpha}(G) = \sum_{uv \in E(G)} (d(u)+d(v))^{\alpha}$ and $\prod_1^*(G) = \prod_{uv \in E(G)} (d(u)+d(v))$, we have
$$\begin{array}{rcl}
&\chi_{\alpha}(CS_n)& = \left\{ \begin{array}{rcl}
6^{\alpha +1,}~~~~~~~~~~~~~~~~~~~~~~~~~~~~~~~~~~~~~~~~~~~~~~~~~ & \mbox{if} &  n =1, \\
2^{\alpha} 3^{\alpha}(n+4)+  3^{2\alpha} (4n-2)+2^{2\alpha}3^{\alpha}(n-2),  & \mbox{if} &  n \geq 2,
\end{array}\right. \\
&\prod^*_1(CS_n)& =  \left\{ \begin{array}{rcl}
46656,~~~~ ~~~~~~~~~~~~~~~~~~~~~~~~~~~~ & \mbox{if} &  n =1, \\
2^{3n}3^{10n-2},~~~~~~~~~~~~~~~~~~~~~~~~~~~& \mbox{if} & n \geq 2.
\end{array}\right.
\end{array}$$

This completes the proof. $\hfill\Box$


\subsection{Hexagonal networks}
It is known that there exist three regular plane tilings with
composition of same kind of regular polygons such as triangular,
hexagonal and square. In the construction of hexagonal networks,
triangular tiling is   used. A hexagonal network of dimension
n is usually denoted as $HX_n$, where $n$ is the number of
vertices on each side of hexagon. The number of vertices in
hexagonal networks $HX_n$  are $3n^2-3n+ 1$ and number
of edges are $9n^2-15n+6.$

\begin{figure}[htbp]
    \centering
    \includegraphics[width=2.5in]{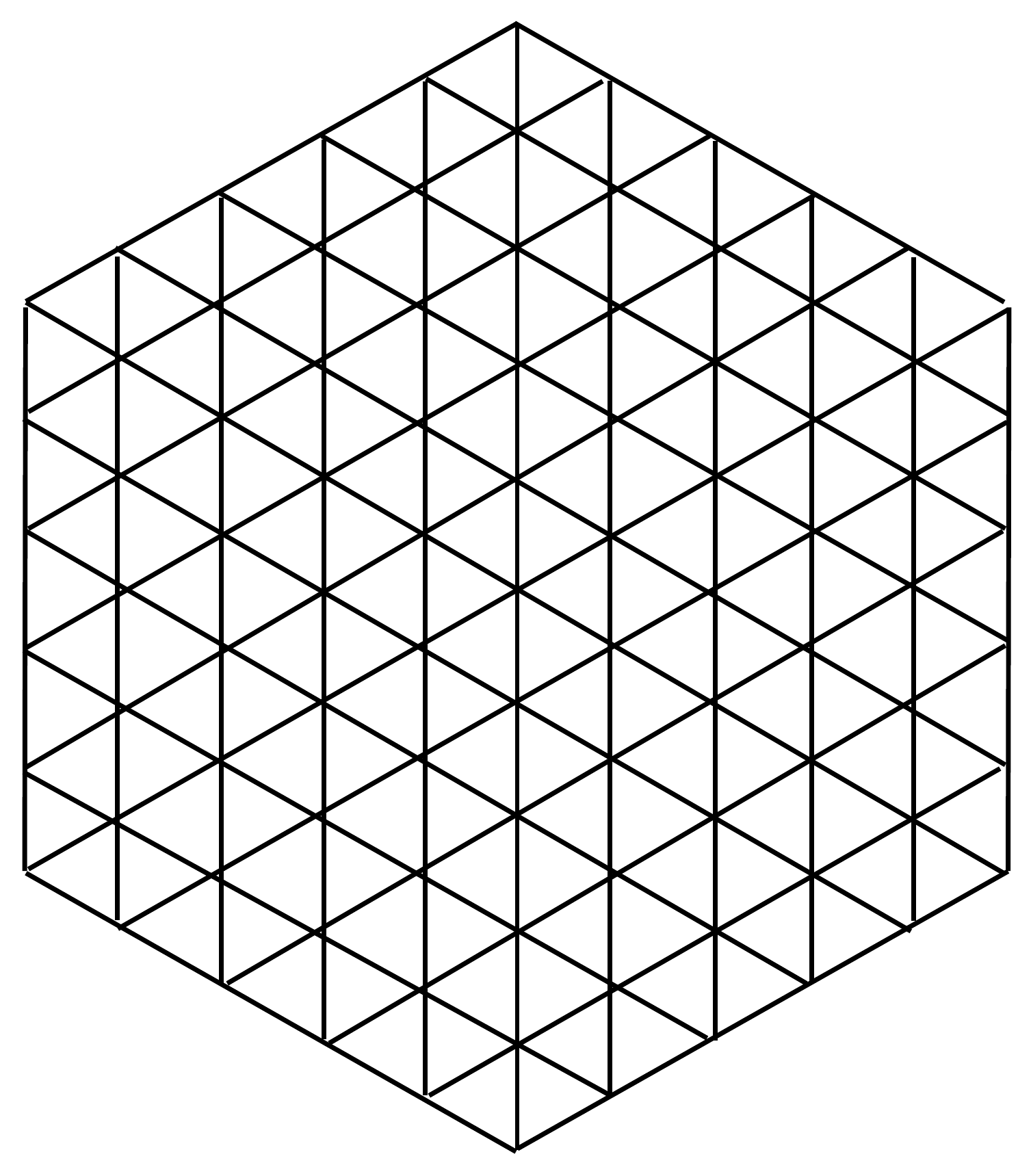}
    \caption{Hexagonal network of dimension 6.}
    \label{fig:sample_figure}
\end{figure}

Here, we compute the exact formulas of
$\chi_{\alpha}(HX_n),$ $\prod_1^*(HX_n),$   $\prod_{1,c}(HX_n)$
and $\prod_2(HX_n)$ for hexagonal networks.

\vskip 2mm \noindent {\bf Theorem 2.3.1. } \emph{Consider the
hexagonal networks $HX_n$, then  the indices of
$\chi_{\alpha}(HX_n),$ $\prod_1^*(HX_n),$ $ \prod_{1,c}(HX_n)$ and
$\prod_2(HX_n)$ are equal to
$$\begin{array}{rcl}
&~\prod_{1,c}(HX_n)&  = 2^{3n^2c+ 3nc -17c}3^{3n^2c-9nc+13c}, \\
& \prod_2(HX_n)& =  2^{18n^2-6n-56}3^{18n^2-54n+60},\\
 &\chi_{\alpha}(HX_n)& = 2^{3\alpha+1}3(n-3) +2^{\alpha+2}3^{\alpha+1}(n-2) - 2^{2\alpha}3^{\alpha+1}(3n^2-11n+10)+7^{\alpha}12+3^{2\alpha}6, \\
&\prod^*_1(HX_n)& = 2^{18n^2-46n-18}3^{9n^2-33n+42}5^{12n-24}7^{12}.
\end{array}$$
}

 \vskip 2mm \noindent  {\bf Proof. }
Let $G$ be the graph of hexagonal network $HX_n$ with $| V(HX_n) |
= 3n^2-3n+ 1$ and $| E(CS_n) | = 9n^2-15n+6$. In a hexagonal
network  $HX_n$,  there are three types of degree of vertices. The
vertices on the corners are of degree 3, the vertices on the
boundary except the corners, are of degree 4 and others are of
degree 6. By Table 5, we get $$n_3=6,~n_4=6n-12,~n_6=3n^2-9n+7.$$
Note that $m_{3,6}=6, m_{3,4} + m_{3,6} +m_{4,4}+m_{4,6}+ m_{6,6}
= 9n^2-15n+6 =\mid E(HX_n)\mid.$ In addition, by Lemma 2.0.1, one
can obtain
$$\begin{array}{rcl}
 m_{3,4} + m_{3,6} &= 3n_3,\\
 m_{3,4} +2m_{4,4} + m_{4,6} &= 4n_4,\\
  m_{3,6}+ m_{4,6} + 2m_{6,6} &= 6n_6.\\
\end{array}$$
Consequently,
$$\begin{array}{rcl}
 m_{3,4} &=& 12,\\
 m_{3,6} &=& 6,\\
 m_{4,4} &=& 6n-18,\\
 m_{4,6} &=& 12n-24,\\
 m_{6,6} &=& 9n^2-33n+30.\\
\end{array}$$
 Table 5 provides the vertex partition of
hexagonal networks $HX_n$.

\begin{table}[thb]
\center
\caption{Vertex partition of $HX_n$ based on degrees of vertices.}
\label{tab:CS}
\vspace{0.2in}
\begin{tabular}{|lc  | c|c|c||c ||}
\hline
 \multicolumn{2}{|c|}{Degrees} &\multicolumn{ 1}{c|}{$\;\;\;\;\;$3$\;\;\;\;\;$}  &\multicolumn{ 1}{c|}{$\;\;\;\;\;$4$\;\;\;\;\;$} &\multicolumn{ 1}{c|}{$\;\;\;\;\;$6$\;\;\;\;\;$} \\
\hline\hline
  & Number of vertices     & \multicolumn{ 1}{c|} {$6$} & \multicolumn{ 1}{c|}{$6n-12$} & \multicolumn{ 1}{c|}{$3n^2-9n+7$}  \\
\hline
\end{tabular}
\end{table}

By Table 5, $\prod_{1,c}(G) = \prod_{v \in V(G)}d(v)^c$ and  $\prod_2(G) = \prod_{u \in V(G)}
d(u)^{d(u)}$, we have

$$\begin{array}{rcl}
&~\prod_{1,c}(HX_n)&  = 2^{3n^2c+ 3nc -17c}3^{3n^2c-9nc+13c}, \\
& \prod_2(HX_n)& =  2^{18n^2-6n-56}3^{18n^2-54n+60}.
\end{array}$$

In a hexagonal network  $HX_n$,  there are five types of edges
based on the degree of the vertices of each edge.  The graphic
properties on the degree sum of vertices for each edge imply the
edge partition of silicate networks $HX_n$. By Lemma 2.0.1.,
Table 6
 explains such partition.

\begin{table}[thb]
\center
\caption{Edge partition of $HX_{n}$ based on degrees of end vertices of each edge.}
\label{tab:CS}
\vspace{0.2in}
\begin{tabular}{|lc |c | c|c|c||c ||}
\hline
 \multicolumn{2}{|c|}{($d(u),d(v)$) where $u,v \in E(G)$} &\multicolumn{ 1}{c|}{$\;\;\;$(3,4)$\;\;\;$}
  &\multicolumn{ 1}{c|}{$\;\;\;$(3,6)$\;\;\;$}  &\multicolumn{ 1}{c|}{$\;\;\;$(4,4)$\;\;\;$}
   &\multicolumn{ 1}{c|}{$\;\;\;$(4,6)$\;\;\;$} & \multicolumn{ 1}{c|}{(6,6)} \\
\hline\hline
  & Number of edges     & \multicolumn{ 1}{c|} {$12$} & \multicolumn{ 1}{c|}{$6$} &
   \multicolumn{ 1}{c|}{$6n-18$} & \multicolumn{ 1}{c|}{$12n-24$}& \multicolumn{ 1}{c|}{$9n^2- 33n +30$} \\
\hline
\end{tabular}
\end{table}

By Table 6, $\chi_{\alpha}(G) = \sum_{uv \in E(G)} (d(u)+d(v))^{\alpha}$ and $\prod_1^*(G) = \prod_{uv \in E(G)} (d(u)+d(v))$, we have
$$\begin{array}{rcl}
  &\chi_{\alpha}(HX_n)& = 2^{3\alpha+1}3(n-3) +2^{\alpha+2}3^{\alpha+1}(n-2) - 2^{2\alpha}3^{\alpha+1}(3n^2-11n+10)+7^{\alpha}12+3^{2\alpha}6, \\
&\prod^*_1(HX_n)& = 2^{18n^2-46n-18}3^{9n^2-33n+42}5^{12n-24}7^{12}.
\end{array}$$

This completes the proof. $\hfill\Box$


\subsection{ Oxide networks}
Oxide networks are of vital importance  in the study of silicate networks. If we delete silicon vertices from a silicate network, we
obtain an oxide network. An $n$-dimensional oxide network is denoted as $OX_n$. The number of vertices  is $9n^2+ 3n$ and the edge set cardinality is $18n^2$.

The exact formulas of
$\chi_{\alpha}(OX_n),$ $\prod_1^*(OX_n),$   $\prod_{1,c}(OX_n)$
and $\prod_2(OX_n)$ for oxide  networks are computed below.

\begin{figure}[htbp]
    \centering
    \includegraphics[width=2.5in]{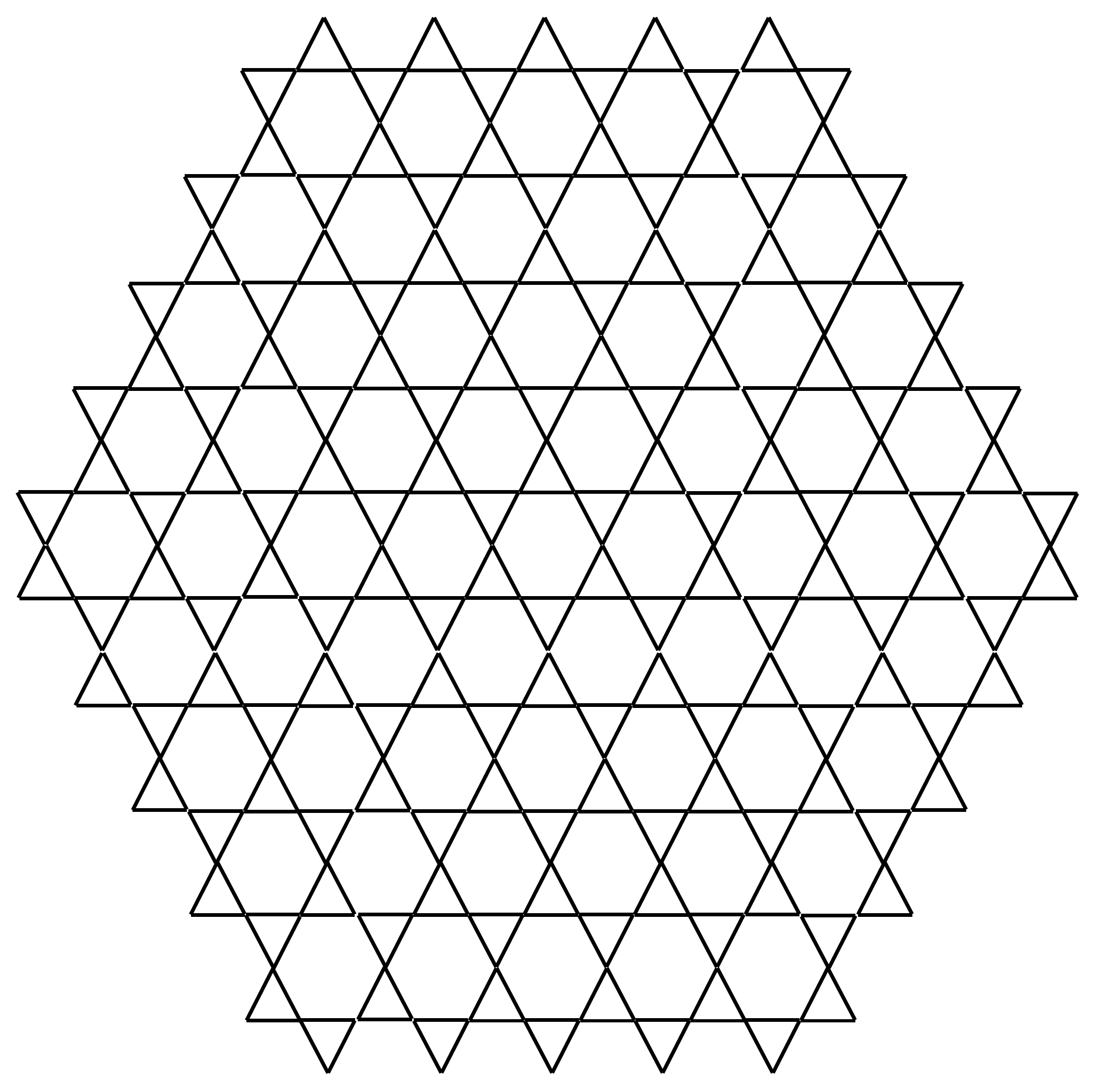}
    \caption{Oxide network of dimension 5.}
    \label{fig:sample_figure}
\end{figure}

\vskip 2mm \noindent  {\bf Theorem 2.4.1. } \emph{Consider the
oxide networks $OX_n$, then  the indices of
$\chi_{\alpha}(OX_n),$ $\prod_1^*(OX_n),$  $\prod_{1,c}(OX_n)$ and
$\prod_2(OX_n)$ are equal to
$$\begin{array}{rcl}
&~\prod_{1,c}(OX_n)&  = 2^{18n^2c}, \\
& \prod_2(OX_n)& =  2^{72n^2-12n},\\
 &\chi_{\alpha}(OX_n)& = n2^{\alpha+2}3^{\alpha+1}+3n(3n-2)2^{3\alpha+1}, \\
&\prod^*_1(OX_n)& = 2^{54n^2-24n}3^{12n}.\\
\end{array}$$
}

 \vskip 2mm \noindent {\bf Proof. }
Let $G$ be the graph of oxide network $OX_n$ with $|V(OX_n) | = 9n^2+
3n$ and $| E(OX_n) | = 18n^2$. Because of the graphic properties on
the degree  of vertices, we have the vertex partition of hexagonal
networks $OX_n$. Table 7  explains such partition for $OX_n$.

\begin{table}[thb]
\center
\caption{Vertex partition of $OX_n$ based on degrees of vertices.}
\label{tab:CS}
\vspace{0.2in}
\begin{tabular}{|lc  | c|c|c||c ||}
\hline
 \multicolumn{2}{|c|}{Degrees} &\multicolumn{ 1}{c|}{$\;\;\;\;\;$2$\;\;\;\;\;$}  &\multicolumn{ 1}{c|}{$\;\;\;\;\;$4$\;\;\;\;\;$}  \\
\hline\hline
  & Number of vertices     & \multicolumn{ 1}{c|} {$6n$} & \multicolumn{ 1}{c|}{$9n^2 - 3n$}  \\
\hline
\end{tabular}
\end{table}

By Table 7, $\prod_{1,c}(G) = \prod_{v \in V(G)}d(v)^c$ and  $\prod_2(G) = \prod_{u \in V(G)}
d(u)^{d(u)}$, we have

$$\begin{array}{rcl}
&~\prod_{1,c}(OX_n)&  = 2^{6cn}4^{c(9n^2-3n)}= 2^{18n^2c},  \\
& \prod_2(OX_n)& =  2^{12n}4^{4(9n^2-3n)} =  2^{72n^2-12n}.
\end{array}$$

Based on the graphic properties and Table 7, we get
$$n_2=6n,~n_4=9n^2-3n.$$
Note that $m_{2,4} + m_{4,4}= 18n^2=\mid E(OX_n)\mid.$ In
addition, by Lemma 2.0.1, one can obtain
$$  m_{2,4}  = 2n_2,~
 m_{2,4} +2m_{4,4} = 4n_4. $$
Consequently,
$$  m_{2,4} = 12n,~
 m_{4,4} = 18n^2-12n. $$
By the above calculation, the edge partition of silicate networks
$OX_n$ are obtained in Table 8.

\begin{table}[thb]
\center
\caption{Edge partition of $OX_n$ based on degrees of end vertices of each edge.}
\label{tab:CS}
\vspace{0.2in}
\begin{tabular}{|lc  | c|c|c||c ||}
\hline
 \multicolumn{2}{|c|}{($d(u),d(v)$) where $u,v \in E(G)$} &\multicolumn{ 1}{c|}{$\;\;\;\;\;$(2,4)$\;\;\;\;\;$}  &\multicolumn{ 1}{c|}{(4,4)}  \\
\hline\hline
  & Number of edges     & \multicolumn{ 1}{c|} {$12n$} & \multicolumn{ 1}{c|}{$18n^2 - 12n$}  \\
\hline
\end{tabular}
\end{table}

By Table 8, $\chi_{\alpha}(G) = \sum_{uv \in E(G)} (d(u)+d(v))^{\alpha}$ and $\prod_1^*(G) = \prod_{uv \in E(G)} (d (u)+d (v))$, we have
$$\begin{array}{rcl}
 &\chi_{\alpha}(OX_n)& = 12n6^{\alpha}+(18n^2-12n)8^{\alpha}= n2^{\alpha+2}3^{\alpha+1}+3n(3n-2)2^{3\alpha+1}, \\
&\prod^*_1(OX_n)& = 6^{12n}8^{18n^2-12n}=2^{54n^2-24n}3^{12n}.
\end{array}$$

This completes the proof. $\hfill\Box$


\subsection{ Honeycomb networks}
Honeycomb networks are very useful in computer graphics, cellular phone base stations, image processing and as a representation
of benzenoid hydrocarbons in chemistry. If we recursively use hexagonal tiling in a particular pattern, honeycomb
networks are formed. An $n$-dimensional honeycomb network is denoted as $HC_n$, where $n$ is the number of hexagons
between central and boundary hexagon. Honeycomb network $HC_n$ is constructed from $HC_{n-1}$ by adding a layer of hexagons
around boundary of $HC_{n-1}$.
The number of vertices are  $6n^2$ and  the number of edges are $9n^2-3n$.

\begin{figure}[htbp]
    \centering
    \includegraphics[width=2.5in]{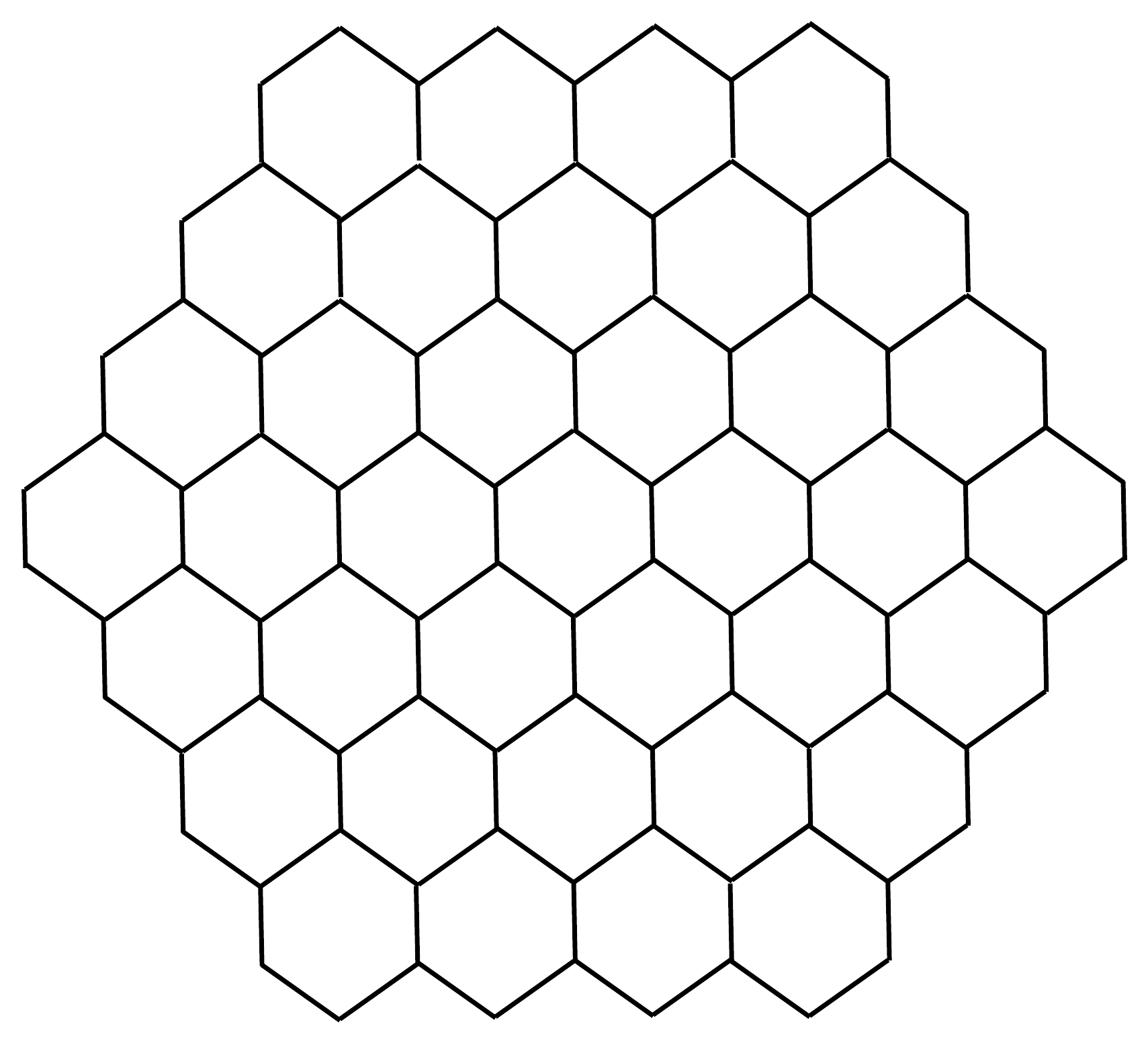}
    \caption{Honeycomb network of dimension 4.}
    \label{fig:sample_figure}
\end{figure}

In the following result,  the exact formulas of
$\chi_{\alpha}(HC_n),$ $\prod_1^*(HC_n),$   $\prod_{1,c}(HC_n)$
and $\prod_2(HC_n)$ for honeycomb networks are calculated.

\vskip 2mm \noindent {\bf Theorem 2.5.1. } \emph{Consider the
honeycomb networks $HC_n$, then  the indices of
$\chi_{\alpha}(HC_n),$ $\prod_1^*(HC_n),$  $\prod_{1,c}(HC_n)$ and
$\prod_2(HC_n)$ are equal to
$$\begin{array}{rcl}
&~\prod_{1,c}(HC_n)&  = 2^{6nc}3^{6n^2c-6nc}, \\
& \prod_2(HC_n)& =  2^{12n}3^{18n^2-18n},\\
 &\chi_{\alpha}(HC_n)~& = 3 \cdot 2^{2\alpha+1}+12(n-1)5^{\alpha}+2^{\alpha}3^{\alpha+1}(3n^2-5n+2), \\
&\prod^*_1(HC_n)& =  2^{9n^2-15n+18}3^{9n^2-15n+6}5^{12(n-1)}.
\end{array}$$
}

 \vskip 2mm \noindent {\bf Proof. }
Let $G$ be the graph of honeycomb network $HC_n$ with $| V(HC_n) | =6n^2$
and $| E(HC_n) |  = 9n^2-3n$. The
vertex partition of honeycomb networks $HC_n$ based on the degree
of vertices are vertices of degree 2 and 3.
The vertices of the boundary in one hexagon are either of degree 2 and others are of degree 3. Table 9 shows such partition for $HC_n$ below.

\begin{table}[thb]
\center
\caption{Vertex partition of $HC_n$ based on degrees of vertices.}
\label{tab:CS}
\vspace{0.2in}
\begin{tabular}{|lc  | c|c|c||c ||}
\hline
 \multicolumn{2}{|c|}{Degrees} &\multicolumn{ 1}{c|}{$\;\;\;\;\;$2$\;\;\;\;\;$}  &\multicolumn{ 1}{c|}{$\;\;\;\;\;$3$\;\;\;\;\;$}  \\
\hline\hline
  & Number of vertices     & \multicolumn{ 1}{c|} {$6n$} & \multicolumn{ 1}{c|}{$6n^2 - 6n$}  \\
\hline
\end{tabular}
\end{table}

By Table 9, $\prod_{1,c}(G) = \prod_{v \in V(G)}d(v)^c$ and  $\prod_2(G) = \prod_{u \in V(G)}
d(u)^{d(u)}$, we have
$$\begin{array}{rcl}
&~\prod_{1,c}(HC_n)&  = 2^{6nc}3^{6n^2c-6nc}, \\
& \prod_2(HC_n)& =  2^{12n}3^{18n^2-18n}.
\end{array}$$

In a honeycomb network  $HC_n$, there are three types of edges
based on the degree of the vertices of each edge. There are three
types of edges of honeycomb network $HC_{n}$ based on degrees of
vertices. Based on the graphic properties and Table 9, we get
$$n_2=6n,~n_3=6n^2-6n.$$
Note that $m_{2,2} + m_{2,3}+ m_{3,3}= 9n^2-3n=\mid E(HC_n)\mid.$
In addition, by Lemma 2.0.1, one can obtain
$$2m_{2,2}+m_{2,3}  = 2n_2,~
 m_{2,3} +2m_{3,3} = 3n_3.$$
Consequently,
$$\begin{array}{rcl}
 m_{2,2} &=& 6,\\
 m_{2,3} &=& 12n-12,\\
 m_{3,3} &=& 9n^2-15n+6.\\
\end{array}$$
By the above calculation, Table 10 shows them with corresponding
partite set cardinalities.

\begin{table}[thb]
\center
\caption{Edge partition of $HC_n$ based on degrees of end vertices of each edge.}
\label{tab:CS}
\vspace{0.2in}
\begin{tabular}{|lc  | c|c|c||c ||}
\hline
 \multicolumn{2}{|c|}{($d(u),d(v)$) where $u,v \in E(G)$} &\multicolumn{ 1}{c|}{$\;\;\;\;\;$(2,2)$\;\;\;\;\;$}  &\multicolumn{ 1}{c|}{(2,3)}  &\multicolumn{ 1}{c|}{(3,3)}\\
\hline\hline
  & Number of edges     & \multicolumn{ 1}{c|} {6} & \multicolumn{ 1}{c|}{$12n-12$} & \multicolumn{ 1}{c|}{$9n^2-15n+6$}  \\
\hline
\end{tabular}
\end{table}

By Table 10, $\chi_{\alpha}(G) = \sum_{uv \in E(G)} (d(u)+d(v))^{\alpha}$ and $\prod_1^*(G) = \prod_{uv \in E(G)} (d(u)+d(v))$, we have
$$\begin{array}{rcl}
 &\chi_{\alpha}(HC_n)~& = 3 \cdot 2^{2\alpha+1}+12(n-1)5^{\alpha}+2^{\alpha}3^{\alpha+1}(3n^2-5n+2), \\
&\prod^*_1(HC_n)& = 4^6 5^{12(n-1)} 6^{9n^2-15n+6} =   2^{9n^2-15n+18}3^{9n^2-15n+6}5^{12(n-1)}.
\end{array}$$

This completes the proof. $\hfill\Box$


\section{Conclusion}
This paper is devoted to the study of certain degree based topological indices of certain chemical networks having both chemical and structural significance.
More precisely, we have considered the silicate networks, hexagonal networks, honeycomb networks and oxide networks for their structural study.  We determined the sum-connectivity index $\chi_{\alpha}(G),$ the multiplicative version of ordinary first
Zagreb index $\prod_1^*(G),$   the first  multiplicative
Zagreb indices $\prod_{1,c}(G)$
and  the  second multiplicative
Zagreb indices $\prod_2(G)$ for these chemical networks.

Furthermore, we provide the direct relation of $\chi_{\alpha}(G),$ $\prod_1^*(G),$
$\prod_{1,c}(G)$ and $\prod_2(G)$ under considered networks.  It
will be quite helpful to understand their underlying topologies. For the sake of refinement of the image, we fix $\alpha$ and $c$ equal to $2$, $n \geq 2$. In the  figures 8-11, the comparisons of these networks with all those indices  are given.

\section{Acknowledgment}
The authors would like to express their sincere gratitude to Dr. Emeric Deutsch and his insightful suggestions, which led to a number of improvements.

\begin{figure}[htbp]
    \centering
    \includegraphics[width=4in]{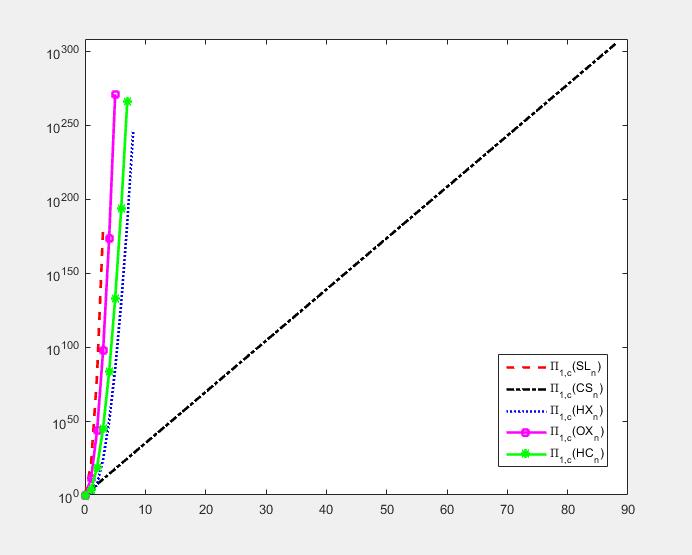}
    \caption{ The computer-based comparative graphs of the first multiplicative Zagreb index for all network families.}
    \label{fig: te}
\end{figure}

\begin{figure}[htbp]
    \centering
    \includegraphics[width=4in]{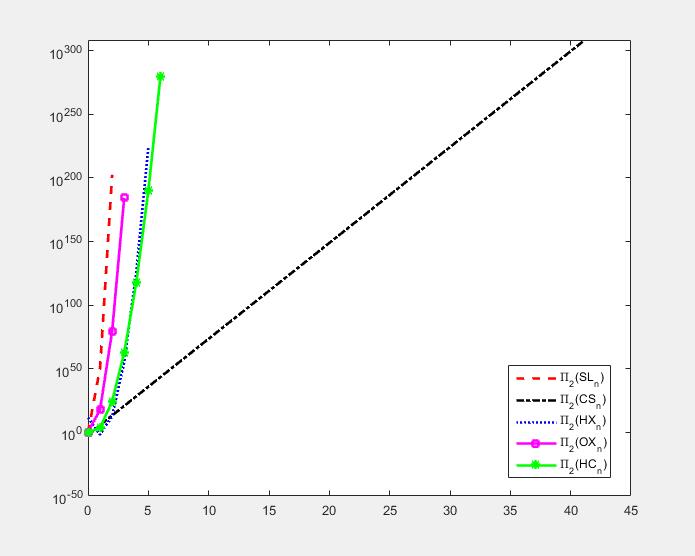}
    \caption{ The comparative graphs of the second multiplicative Zagreb index for all the families of networks.}
    \label{fig: te}
\end{figure}

\begin{figure}[htbp]
    \centering
    \includegraphics[width=4in]{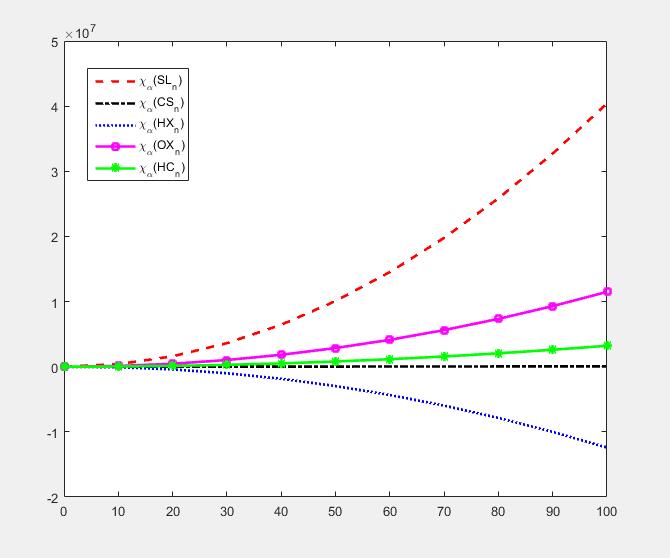}
    \caption{ The comparisons for the sum-connectivity index for all network families.}
    \label{fig: te}
\end{figure}

\begin{figure}[htbp]
    \centering
    \includegraphics[width=4in]{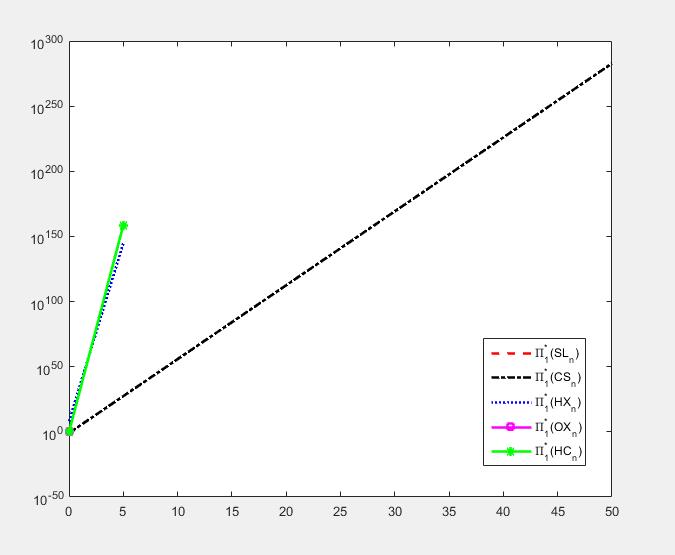}
    \caption{ The graphical representations of multiplicative version of ordinary first Zagreb index for all network families.}
    \label{fig: te}
\end{figure}


\end{document}